\documentclass[11pt,psfig]{article}
\usepackage{amssymb,amsmath,amscd}
\usepackage[all]{xy}
\usepackage{mathtools}
\newtheorem{lem}{Lemma}[section]
\newtheorem{thm}{Theorem}[section]
\newtheorem{cor}{Corollary}[section]

\newcommand{\f}[1]{\mathfrak{#1}}

\newcommand{\Ind}{{\rm Ind}}

\newcommand{\commentout}[1]{}

\newcommand{\mc}{\mathcal}

\newcommand{\diag}{{\,\rm diag}}
\newcommand{\Hom}{{\rm {Hom}\,}}
\newcommand{\sgn}{\, {\rm sgn}}

\newcommand{\vertiii}[1]{{\left\vert\kern-0.25ex\left\vert\kern-0.25ex\left\vert #1 
    \right\vert\kern-0.25ex\right\vert\kern-0.25ex\right\vert}}

\begin{document}
\title{Certain $L^2$-norm and Asymptotic bounds of Whittaker functions for $GL(n, \mathbb R)$}
\author{Hongyu He  \\
Department of Mathematics \\
Louisiana State University \\
email: hongyu@math.lsu.edu\\
}
\date{}
\maketitle
\commentout{\footnote{Key word:  Whittaker model, Jacquet Whittaker function, Kirilov model, $GL(n)$, Asymptote, cusp forms, Harish-Chandra module, unitary representation, perturbation of Hilbert structure}}
\abstract{Whittaker functions of $GL(n, \mathbb R)$ (\cite{jw}\cite{ja}), are most known for its role in the Fourier-Whittaker expansion of cusp forms ( \cite{sh} \cite{ps}). Their behavior in the Siegel set, in large, is well-understood. In this paper, we insert into the literature some potentially useful properties of Whittaker function over the group $GL(n, \mathbb R)$ and the mirobolic group $P_n$. We proved the square integrabilty of the Whittaker functions with respect to certain measures, extending a theorem of Jacquet and Shalika (\cite{js}). For principal series representations, we gave various asymptotic bounds of smooth Whittaker functions over the whole group $GL(n, \mathbb R)$. Due to the lack of good terminology, we use whittaker functions to refer to $K$-finite or smooth vectors in the Whittaker model.  }

\section{Introduction}
Let $G=GL(n)=GL(n, \mathbb R)$. Let $\overline{U}$ be the group of unipotent lower triangular matrices, $U$ be the group of unipotent upper triangular matrices, $A$ be the group of positive diagonal matrices, and $M$ be the centralizer of $A$ in $G$. 
Let $K$ be the standard orthogonal group $O(n)$. Then we have the Iwasawa decomposition $G=KA U$. \\
\\
Let $(\pi, \mc H)$ be a irreducible Hilbert representation of $G$. Let $\pi^{\infty}$ be the Frechet space of smooth vectors in $\mc H$. Let $(\pi^{*})^{-\infty}$ be the topological dual space of $\pi^{\infty}$. This dual space is equipped with the natural action of $\pi^*$.
A functional $\psi \in (\pi^*)^{-\infty}$ is called a Whittaker functional if 
 $$\pi^*(u) \psi= \exp - 2 \pi i (\sum m_i u_{i, i+1}) \psi, \qquad ( \forall \,\, u \in U)$$ for some $m \in \mathbb R^{n-1}$ with $\prod m_i \neq 0$ (\cite{jw} \cite{ja}). In \cite{sh}, J. Shalika showed that Whittaker functionals, if exist, are unique up to a constant for $\pi^{\infty}$. If a Whittaker functional exists, the representation $\pi$ is often said to be generic. In \cite{ko}, Kostant showed that a representation is generic if and only if it has the maximal Gelfand-Kirillov dimension (\cite{ko}).
  We then know the classification of generic representations for $GL(n)$, namely, those irreducible representations induced from the discrete series of $GL(1)$ and $GL(2)$ Levi factors.\\
 \\
  Now fix an $m$. We can define
 $$Wh_{, m}: f \in \pi^{\infty} \rightarrow Wh{,m}(f)=\langle \pi(g) f, \psi_m \rangle.$$
 If $m=\mathbf 1$, the constant vector with entries $1$, we write
 $$Wh_f(g)= \langle \pi(g) f, \psi_{\mathbf 1} \rangle.$$  
 Obviously, $Wh_f \in C^{\infty}( \mathbb C_{\mathbf 1} \times_{U} G)$. Now the group $G$ will act from the right. The space $\{ Wh_f \}$ is known as the Whittaker model of $\pi$, introduced by Jacquet in his studies of automorphic forms (\cite{ja}). Because of the Iwasawa decomposition, $Wh_f(g)$ can be uniquely determined by its restriction on $AK$. Hence Whittaker model is often regarded as a space of smooth functions on $AK$. If $f$ is a spherical vector, then $Wh_f(g)$ is uniquely determined by $Wh_f|_A$. Generally speaking, 
 $ Wh_{f, m}(g) $ is a smooth section of 
 $$\mathbb C_ m \times_{U} G \rightarrow  U \backslash G$$
 where $\mathbb C_{m}$ denote the one dimensional representation defined by the character of $U$:
 $$u \rightarrow \exp 2 \pi i \sum_{i=1}^{n-1} u_{i, i+1} m_i.$$
\noindent 
 Let $P_n$ be the mirabolic subgroup consisting of invertible  matrices with last row $(0,0, \ldots, 0, 1)$. $P_n$ is the semidirect product of $GL(n-1)$ and $\mathbb R^{n-1}$. When we refer to $GL(n-1)$ as a subgroup of $GL(n)$, $GL(n-1)$ will always lie in the upper left corner. Consider the restriction $Wh_{f}|_{P_n}$. It is a smooth section of 
  $$\mathbb C_{\mathbf 1} \times_{U} P_n \rightarrow  U \backslash P_n,$$
  which can be identified with a smooth section of 
  $$\mathbb C_{\mathbf 1} \times_{U_{n-1}} GL(n-1) \rightarrow  U_{n-1} \backslash GL(n-1)$$
  and vice versa. Fix an invariant measure on $GL(n-1)$ and  the Euclidean measure on $\mathbb R^{n-1}$, and  equip $P_n$ with the right invariant measure. Jacquet and Shilika proved the following theorem.
  \begin{thm}[  Jacquet-Shalika \cite{js}]
Let $(\pi, \mc H)$ be a generic irreducible unitary representation of $G$. Then  $$Wh_f|_{P_n} \in L^2(\mathbb C_{\mathbf 1} \times_{U} P_n).$$
The Whittaker model produces a unitary equivalence between
$(\pi|_{P_{n}}, \mc H)$ and $(R, L^2(\mathbb C_{\mathbf 1} \times_{U} P_n))$. Here  $R$ stand for the right regular representation.
\end{thm}
By Mackey's theory, $(R, L^2(\mathbb C_{\mathbf 1} \times_{U} P_n))$ is an irreducible unitary representation of $P_n$, essentially, the unique one with maximal Gelfand-Kirillov dimension (\cite{be} \cite{sa}). This implies that $(\pi|_{P_n}, \mc H)$ is already irreducible, which was conjectured by Kirillov to be true for all irreducible unitary representations of $G$. In the literature, $\{ Wh_f |_{P_n} \}$ is often known as the Kirillov model. \\
\\
Here are our main results.
\begin{thm}[A] Let $\pi$ be an irreducible unitary representation of $GL(n)$ with a Whittaker model. Let $f$ be a $K$-finite vector in $\pi^{\infty}$. Write $a=\diag(a_1, a_2, \ldots, a_{n-1}, 1)$ in the $UAK$ decomposition of $GL(n-1)$. Then for any $l_i \in \mathbb N, i \in [1, n-1]$, we have
$$\prod_{i=1}^{n-1} (\frac{a_i}{a_{i+1}})^{l_i} Wh_f|_{P_n} \in 
L^2(\mathbb C_{\mathbf 1} \times_{U} P_n).$$
Consequently
$$\prod_{i=1}^{n-1} (a_i)^{l_i} Wh_f|_{P_n} \in 
L^2(\mathbb C_{\mathbf 1} \times_{U} P_n).$$
\end{thm}
Identify $L^2(\mathbb C_{\mathbf 1} \times_{U} P_n)$ with 
$L^2(\mathbb C_{\mathbf 1} \times_{U_{n-1}} GL(n-1))$. Here $U_{n-1}= U \cap GL(n-1)$. We have
\begin{thm}[B]
Let $\pi$ be an irreducible unitary representation of $GL(n)$ with a Whittaker model. Let $f$ be a $K$-finite vector in $\pi^{\infty}$. Then
$$Wh_f|_{GL(n-1)} \in L^2(\mathbb C_{\mathbf 1} \times_{U_{n-1}} GL(n-1), |\det|^{\epsilon} d g) $$
for any $\epsilon \geq 0$.
\end{thm}
Indeed, $L^2(\mathbb C_{\mathbf 1} \times_{U_{n-1}} GL(n-1), |\det|^{\epsilon} d g) $ provides us a unitary structure of the perturbed representation $\pi \otimes |\det|^{\frac{\epsilon}{2}}|_{P_n}$. Notice that $\pi \otimes |\det|^{\frac{\epsilon}{2}}$ is never unitary unless $\epsilon \in i \mathbb R$. \\
\\
Theorem A and B suggest that the asymptotic behavior of the Whittaker function can be similarly understood. The asymptotic behavior of the Whittaker function $Wh_f(g)$ is well-know for $g$ in the Siegel set. However, $wh_f(g)$ outside the Siegel set can also be important. In this paper, we prove the following
\begin{thm}[C] Let $\pi(v,\sigma)$ be the principal series defined over $G/MA\overline U$. Suppose that
$$\Re(v_1) < \Re(v_2) \ldots < \Re(v_{n}).$$
Let $\rho=(\frac{n-1}{2}, \frac{n-3}{2}, \ldots, -\frac{n-3}{2}, -\frac{n-1}{2})$ be the half sum of positive roots for $(\f a, \f u)$.
For any $f \in \pi(v, \sigma)^{\infty}$, any $l_i \in \mathbb N, i \in [1, n-1]$,  the function
$$a^{-\rho-v} \prod_{i=1}^{n-1} 
(\frac{a_i}{a_{i+1}})^{l_i} Wh_f(g)$$
is bounded by a constant dependent on $f$ and $l$.
\end{thm}
Hence $$ Wh_f(g) \leq C_{f, l} a^{\rho+v} \prod_{i=1}^{n-1} 
(\frac{a_{i+1}}{a_{i}})^{l_i}.$$
The restriction imposed on $v$, simply means that $v$ is in the open negative Weyl chamber of $\f a^*_{\mathbb C}$.
We shall remark that the asymptotic bound we obtain is the best possible bound. Notice that if $\sigma$ is trivial and $f$ is spherical, as $\frac{a_i}{a_{i+1}} \rightarrow 0$ for all possible $i$, the Whittaker function $Wh_f(a)$ will approach $c(v) a^{\rho+v}$ with $c(v)$ Harish-Chandra's $c$ function (\cite{hel}, \cite{knapp}).  It can be easily verified that $c(v)$ is not zero when $v$ is in the open negative Weyl chamber. Hence $a^{\rho+v}$ is the best possible exponent. Over the Siegel set with $\frac{a_i}{a_{i+1}} \rightarrow \infty$, our theorem implies that $Wh_f(g)$ is fast decaying. More generally, 
 Theorem C  gives effective bounds over  $g=nak$ with $a$ in all Weyl chambers. \\
 \\ 
 At the boundary where $\Re(v_i)=\Re(v_{i+1})$ for some $i$, our theorem is no longer true in general. But a weaker version (with introduction of certain small $\delta$, or $\log$-terms) is true. In the general situation where $\pi$ is induced from  discrete series of $GL(2)$ and $GL(1)$, one can  embed $\pi$ into a principal series $\pi(v, \sigma)$ with $v$ in the closed negative Weyl chamber. Similar statement remains true. But it is unlikely that we can obtain the best results this way. In fact, one would have to use the leading exponents of $\pi$ to state Theorem C correctly (\cite{knapp}, \cite{wa}). There seems to be a deep connection between the asymptotes of the Whittaker functions and leading exponents of $\pi$. We shall not pursue the asymptotes of the more general Whittaker functions involving the discrete series of $GL(2)$ in this paper.\\
\\
Finally, we shall  remark that Theorem C holds for principal series $\pi(v, \sigma)$ of all semisimple Lie groups. 
 \begin{thm}[D]
 Let $G$ be a semisimple Lie group and $NAK$ its Iwasawa decomposition. Let $\overline N$ be the opposite nilpotent group. Let $\Delta^+(\f g, \f a)$ be the positive restricted simple roots. Let $\pi(v, \sigma)$ be a principal series representation built on $G/MA \overline{N}$. Suppose that $\Re (\alpha_i, v) < 0$ for every $\alpha_i \in \Delta^+(\f g, \f a)$. Then for any $f \in \pi(v, \sigma)^{\infty}$, any $l_i \in \mathbb N, i \in [1, \dim A]$,  the function $Wh_f(g)$ is bounded by a constant multiple of 
$$a^{\rho+v+ \sum_{\alpha_i \in \Delta^+(\f g, \f a)} -l_i \alpha_i }.$$

\end{thm}

Unless otherwise stated, all $GL(n)$ in this paper will refer to $GL(n, \mathbb R)$. For $n=2$, the Whittaker functions of $GL(2)$ are classical whittaker functions and their asymptotes are well-known (\cite{jw}). We shall from now on assume $n \geq 3$.

\section{$L^2$-norms of Whittaker Model}
Let $G=GL(n, \mathbb R)$.

 Fix the invariant measure on $GL(n-1)=U_{n-1} A_{n-1} K_{n-1}$ as
$$a^{-2 \rho_{n-1}} (\prod_{1 \leq i < j \leq n-1} d u_{i,j}) (\prod_{i=1}^{n-1} \frac{ d a_i}{ a_i}) d k=a^{-2 \rho_{n-1}} d u \frac{ d a}{a} d k ,$$
with $\rho_{n-1}=(\frac{n-2}{2}, \frac{n-4}{2}, \ldots, -\frac{n-2}{2})$ and $ d a_i$ the Euclidean measure. This induces the right invariant measure on 
$ U_{n-1} \backslash GL(n-1) \cong A_{n-1} K_{n-1}$:
$$a^{-2 \rho_{n-1}} \frac{ d a}{a} dk.$$
Suppose that $(\pi, \mc H)$ is an irreducible unitary representation of $G$. Let $\pi^{\infty}$ be the irreducible  smooth representation associated with $\pi$.
Let $\psi_m$ be a Whittaker functional in $(\pi{^*})^{-\infty}$.
Let $\mc H_{K}$ be the space of $K$-finite vectors in $\mc H$. Then $\mc H_{K}$ is a Harish-Chandra module.
By the theorem of Jacquet-Shalika, for any $f \in \mc H$, 
$$Wh_{f, m}|_{A_{n-1} K_{n-1}} \in L^2(A_{n-1} K_{n-1}, a^{-2 \rho_{n-1}}   \frac{ d a}{ a} d k ).$$
We prove the following 
\begin{thm}\label{apart}
Let $(\pi, \mc H)$ be an irreducible unitary representation and $\psi_{m}$ be a Whittaker functional. For any $f \in \mc H_{K}$,
$$Wh_{f, m}|_{A_{n-1}} \in L^2(A_{n-1}, a^{-2 \rho_{n-1}}  \frac{ d a}{ a}) .$$
\end{thm}
\noindent
Proof: Our proof is standard and depend on the fact that the $\tau$-isotypic subspace $\mc H_{\tau}$ is finite dimensional for every $\tau \in \hat K$. \\
\\
Notice that $\mc H_{K}= \sum_{\tau \in \hat{K}} \mc H_{\tau}$. Without loss of generality, assume that $ f \in \mc H_{\tau}$ and $f$ is in $V_{\tau}$, {\bf an irreducible representation of $K$}. Choose an orthonormal basis $\{ e_i \}$ in $V_{\tau}$.  Then $\pi(k) f=  \sum (\pi(k)f, e_i) e_i$ and
$$ Wh_{f, m}( ak)= \bigl\langle \sum (\pi(k) f, e_i) \pi(a) e_i, \psi_m \bigr\rangle
=\sum (\pi(k) f, e_i) Wh_{e_i, m}(a).$$
Now $\{ (\pi(k)f, e_i) \}$ are orthogonal to each other in $L^2(K)$, but not necessarily in $L^2(K_{n-1})$. However, $V_{\tau}|_{K_{n-1}}$ is multiplicity free (\cite{gw}). We may choose the orthonormal basis $\{E_i \}$ from each irreducible subrepresentation in $V_{\tau}|_{K_{n-1}}$. Now
$ \{ (\pi(k)f, E_i) \}$ is an orthogonal set in $L^2(K_{n-1)})5$ by Schur's lemma.
We then have
$$\int_{A_{n-1} K_{n-1}} | Wh_{f, m} (a k)|^2 d k a^{-2 \rho_{n-1}}  \frac{ d a}{ a}$$
$$= 
\sum \| (\pi(k) f, E_i) \|^2_{L^2(K_{n-1})} \| Wh_{E_i, m}|_{A_{n-1}} \|^{2}_{L^2(A_{n-1}, a^{- 2 \rho_{n-1}} \frac{d a}{a})} .$$
Hence 
$$Wh_{E_i, m}|_{A_{n-1}} \in L^2(A_{n-1}, a^{-2 \rho_{n-1}}  \frac{d a}{a}) .$$
It follows that 
$$Wh_{f, m}|_{A_{n-1}} \in L^2(A_{n-1}, a^{-2 \rho_{n-1}} \frac{ d a}{a}) .$$
$\Box$\\
\\
Let $\{ X_{i,j} \}$ be the standard basis for the Lie algebra  $\f{gl}(n)$.We have
\begin{lem}\label{diff}
$\pi^*(X_{i,i+1}) \psi_m= 2 \pi i m_{i} \psi_m$
\end{lem}
Proof: We have
$$\pi^*(X_{i, i+1}) \psi_m(u)=\frac{ d}{ d t}|_{t=0} \pi^*(\exp t X_{i,i+1})\psi_m$$
$$=(\frac{ d}{ d t}|_{t=0} \exp 2 \pi i m_i t )\psi_m=2 \pi i  m_i \psi_m(u).$$
$\Box$ \\
\\
{\bf Proof of Theorem A}: 
For every $f \in \pi^{\infty}$, we compute
\begin{equation}
\begin{split}
Wh_{\pi(X_{i, i+1})f,m}(a) = 
& \langle \pi(a) \pi(X_{i,i+1}) f, \psi_m \rangle \\
= & \langle \pi(\frac{a_{i}}{a_{i+1}}X_{i, i+1}) \pi(a) f,  \psi_m \rangle \\
= & \frac{a_{i}}{a_{i+1}} \langle \pi(X_{i, i+1}) \pi(a) f, \psi_{ m} \rangle \\
= & \frac{a_{i}}{a_{i+1}} \langle \pi(a) f, - \pi^*(X_{i, i+1}) \psi_{ m} \rangle \\
= & - \frac{a_{i}}{a_{i+1}} 2 \pi i m_i \langle \pi(a)f, \psi_{ m} \rangle \\
= & - 2 \frac{a_i}{a_{i+1}}  \pi i m_i Wh_{f, m}(a)
\end{split}
\end{equation}
Without loss of generality, assume that $f \in \mc H_{\tau}$ for some $\tau \in \hat K$. Then $\pi(X_{i, i+1}) f$ will also be $K$-finite. By Theorem \ref{apart}, $$\frac{a_i}{a_{i+1}} m_i Wh_{f, m}|_{A_{n-1}} \in L^2(A_{n-1}, a^{-2\rho_{n-1}} \frac{ d a}{a}).$$
By induction for any $l_i \in \mathbb N$,
$$\prod_{i=1}^{n-1} (\frac{a_i}{a_{i+1}})^{l_i} Wh_f|_{A_{n-1}} \in 
L^2(A_{n-1}, a^{-2\rho_{n-1}} \frac{ d a}{a}).$$
Let $\{ e_j \, \, j \in [1, \dim \mc H_{\tau}] \}$ be an orthonormal basis of $\mc H_{\tau}$. The above statement is true for every $e_j$. Then
\begin{equation}
\begin{split}
 & \prod_{i=1}^{n-1} (\frac{a_i}{a_{i+1}})^{l_i}  Wh_f(a k) \\
 = & \prod_{i=1}^{n-1} (\frac{a_i}{a_{i+1}})^{l_i} \langle \pi(a k) f, \psi_{\mathbf 1} \rangle \\
= & \prod_{i=1}^{n-1} (\frac{a_i}{a_{i+1}})^{l_i} \langle \pi(a) \sum_{j} (\pi(k) f, e_j) e_j, \psi_{\mathbf 1} \rangle \\
= & \prod_{i=1}^{n-1} (\frac{a_i}{a_{i+1}})^{l_i} \sum_{j} (\pi(k) f, e_j) \langle \pi(a)  e_j, \psi_{\mathbf 1} \rangle
\end{split}
\end{equation}
Since each $(\pi(k) f, e_j) \in C^{\infty}(K_{n-1}) \subseteq L^2(K_{n-1})$, we have $$\prod_{i=1}^{n-1} (\frac{a_i}{a_{i+1}})^{l_i} Wh_f|_{A_{n-1} K_{n-1}} \in 
L^2(A_{n-1} K_{n-1}, a^{-2\rho_{n-1}} \frac{ d a}{a} d k).$$
Since $  U \backslash P_n \rightarrow U_{n-1} \backslash GL(n-1) \cong A_{n-1} K_{n-1}$, The first statement is Theorem A is proved. We also have 
$$a_i^{l_i} Wh_f= \prod_{j=i}^{n-1} (\frac{a_j}{a_{j+1}})^{l_i} Wh_f \in L^2(\mathbb C_{\mathbf 1} \times_U P_n) \qquad (a_n=1).$$
The second statement of Theorem A follows immediately. $\Box$\\

\section{Asymptotic bound on Whittaker model}
The asymptotes of Whittaker functions over the Siegel set in well-known. Let
$$S(t)= \Gamma_U \backslash U A(t) K$$
 be the Siegel set with $\Gamma_U$ a lattice in $U$, $t >0$ and
 $$A(t)=\{a \in A \mid \frac{a_i}{a_{i+1}} \geq t \qquad (\forall \,\,\, i \in [1, n-1]) \}.$$
 Let $\pi$ be an irreducible smooth representation with a Whittaker model and $f \in \pi^{\infty}$. Then for any $l \in \mathbb N^{n-1}$ and $ u a k \in S(t)$,
 \begin{equation}\label{fast}
 |Wh_{f}(u a k)| \leq C_{l, f} \prod_{i=1}^{n-1} ( \frac{a_i}{a_{i+1}})^{-l_i}.
 \end{equation}
 Notice that, with respect to the right invariant measure, the Siegel set has a finite measure. Hence this estimate is inadequate for the discussion of $Wh_{f}(g)$ with $g \in G$.\\
 \\
 The fast decaying property
 in Equation \ref{fast} is also a special feature  of cusp forms of $G$.
  Indeed, in the main lemma of Ch.1 \cite{hc}, Harish-Chandra proved that if a smooth function of $G$ with  zero constant terms along the directions of maximal parabolics containing $NA$, then this function is fast decaying on the Siegel set. Harish Chandra's intention was to apply this lemma to cusp forms. It is not hard to see that the Whittaker function $Wh_{f}(g)$ has zero constant terms for all parabolics containing $NA$. Hence Equation \ref{fast} also follows from Harish-Chandra's lemma.
  Philosophically, this observation provides us some initial evidence that many analytic properties of the Whittaker functions hold for cusp forms and vice versa.\\
  \\
  We shall also remark that fast decaying property over Siegel set is often rephrased simply by using the norm $\| g \|^{-n}$ with $g \in S$ (\cite{gs}). This is no longer good for the purpose of studying $Wh_f(g)$ for arbitrary $g$. We have to go back to the classical Iwasawa decomposition.
  
  \subsection{General setup of Principal series representation} 
 Fix a minimal parabolic subgroup $MA\overline{U}$, consisting of all invertible lower triangular matrices. 
Form the principal series representation $\pi(v,\sigma)$ with $g \in G$ acting from the left on smooth functions of the form
$$f(g ma \overline{u})= \sigma(m)^{-1} a^{\rho-v} f(g)$$
where
\begin{enumerate}
\item $\pi(v, \sigma)(g) f(x)=f(g^{-1} x)$;
\item $\rho=(\frac{n-1}{2}, \frac{n-3}{2}, \ldots, -\frac{n-1}{2})$ is the half of the positive roots of $(\f u, \f a)$;
\item $v \in \Hom( \f a, \mathbb C)=\mathbb C^n$ is identified with the characters of $A$ and $a^v=\prod_{i=1}^n a_i^{v_i}$;
\item $\sigma=(\sigma_1, \sigma_2, \ldots, \sigma_n)$ is a character of $M=\{ \diag(\pm 1, \pm 1, \ldots \pm 1) \}$.
\end{enumerate}

Each $\pi(v,\sigma)$ restricted to the center of $GL(n)$, produces a central character
$$\diag(a_1,a_1, \ldots, a_1) \rightarrow \prod_{i=1}^n \sigma_i(\sgn(a_1))| a_1|^{v_i}.$$
We may speak of representations of central character $\chi$. 
We use $\pi(v,\sigma)^{\infty}$ to denote the linear space of all smooth vector $f$.\\
\\
The function $f$ is uniquely determined by its value on $K$ and vice versa. Hence the smooth representation $\pi^{\infty}_{v, \sigma}$ can be identified with smooth sections of the vector bundle
$$ K \times_M \mathbb C_{\sigma}.$$
This is often called the (smooth) compact picture, or the compact model. Fix an invariant  measure on $K$. We now equip 
$\pi(v,\sigma)^{\infty} \cong C^{\infty}(K \times_M \mathbb C_{\sigma}, K/M)$ with the $L^2$-norm on $K/M$. Then $\pi(v,\sigma)$ 
 becomes a unitary representation when  $v $ is purely imaginary. We use $( *, *)$ to denote the Hilbert inner product
associated with  $L^2(K \times_M \mathbb C_{\sigma})$. In addition, there is a  canonical (complex linear) non-degenerate pairing between $\pi(v,\sigma)$ and $\pi(-v,\sigma^*)$
$$\langle f_1, f_2 \rangle=\int_{K/M} \langle f_1(k), f_2(k) \rangle d [k].$$
Here $\sigma^* \cong \sigma$ and $ \langle f_1(k), f_2(k) \rangle= f_1(k) f_2(k)$ only depends on $[k] \in K/M$. It follows that $\pi(-v,\sigma^*)$ can be identified with $ \pi(v,\sigma)^*$,
 the dual representation of $\pi(v,\sigma)$. \\
 \\
 Let $\| f \|$ denote the $L^2$-norm $\| f \|_{L^2(K \times_M \mathbb C_{\sigma})}$. 
 We equip $\pi^{\infty}(v, \sigma)$ with the semi-norms $\| f \|_{X}= \| \pi(v,\sigma)(X) f \|$ with $X \in U(\f g)$. Then $\pi(v,\sigma)^{\infty}$ becomes a Frechet space. Its dual space, consisting of continuous linear functionals, contains $\pi(-v,\sigma^*)^{\infty}$ as a subspace. It is often denoted by $\pi(-v,\sigma^*)^{-\infty}$.
 We retain $\pi(v,\sigma)(g) (g \in G)$ and $\pi(v,\sigma)(X) (X \in U(\f g))$ for the group action and Lie algebra action on $\pi(v,\sigma)^{-\infty}$. 
 \subsection{Noncompact model and Proof of Theorem C}
 By the Bruhat decomposition, the image of $U$ in $G/MA\overline{U}$ is open and dense. We consider the pull back 
$$i^*: \pi(v,\sigma)^{\infty} \rightarrow C^{\infty}(U).$$
This is essentially the restriction map of $f$ to the group $U$. It is injective. The group $G$ acts on the image of $i^*$ conformally. $ i^*(\pi(v,\sigma)^{\infty})$ is often called the noncompact model. Fix the invariant measure, the Euclidean measure on $[ u_{ij} ]$. If $v$ is purely imaginary, then $\pi(v,\sigma)(g)$ will be unitary operators on $L^2(U)$. Hence we have a unitary representation $(\pi(v,\sigma), L^2(U))$. Similarly, we have a nondegenerate  pairing between $\pi(v,\sigma)^{\infty}$ and $\pi(-v,\sigma^*)^{\infty}$:
$$\langle f_1, f_2 \rangle=\int_{U} \langle f_1(u), f_2(u) \rangle \prod_{i < j} d u_{i,j} .$$
We normalize the invariant measure on $K$ so the pairing defined here is the same as the pairing defined over $K/M$. Throughout this paper, we shall use the noncompact model unless otherwise stated. \\
\\
The following lemma is well-known and crucial to the construction of the standard intertwining operator and Jacquet's Whittaker function.
\begin{lem}\label{l1noncompact} Suppose that $\Re(v_1) < \Re(v_2) \ldots < \Re(v_n)$. Let $f \in \pi(v, \sigma)^{\infty}$. Then $f_|K$ is bounded and $f|_N \in L^1(N)$.
\end{lem}
Indeed for any $X \in U(\f g)$, $\pi(v, \sigma)(X) f|_K$ will also be bounded. 
See for example Ch. VII. 8,9,10 in \cite{knapp}. 
Write 
$$({\f a}^*_{\mathbb C})^{-}=\{ v \in \f a^*_{\mathbb C} \mid \Re(v_1) < \Re(v_2) \ldots < \Re(v_n) \}.$$
We call this open negative Weyl chamber in $\f a^*_{\mathbb C}$. The closed negative Weyl chamber
$$cl({\f a}^*_{\mathbb C})^{-}=\{ v \in \f a^*_{\mathbb C} \mid \Re(v_1) \leq \Re(v_2) \leq \ldots \leq \Re(v_n) \}.$$

Now we take $\psi_m$ to be the function $ \exp 2 \pi i (\sum_{i=1}^{n-1} m_i u_{i, i+1})$ when $v \in ({\f a}^*_{\mathbb C})^{-}$.   Since all vectors in the noncompact model $i^*(\pi(v, \sigma)^{\infty})$ are in $L^1(U)$,  the function $ \exp 2 \pi i (\sum_{i=1}^{n-1} m_i u_{i, i+1})$ defines a continuous linear functional on $\pi(v, \sigma)^{\infty}$. We then analytically continue $\psi_m$ of $\pi(v, \sigma)$ to all $cl({\f a}^*_{\mathbb C})^{-}$.

\begin{lem} We have
 $\pi(-v, \sigma)(a^{-1}) \psi_m=a^{\rho+v} \psi_{Ad(a)m}$ where 
 $$Ad(a)m=(\frac{a_1}{a_2}m_1, \frac{a_2}{a_3} m_2, \ldots, \frac{a_{n-1}}{a_n}m _{n-1}).$$
\end{lem}
For the reason of bookkeeping, we check that
$$\pi(-v, \sigma)(a^{-1}) \psi_m (u)=\psi_m(a u)=\psi_m( Ad(a)u a)$$
$$=a^{-(-\rho-v)} \psi_m(Ad(a) u)= a^{\rho+v} \psi_{Ad(a) m}(u).$$
$\Box$ \\
\\
 Let $f_0$ be the spherical vector of $\pi(\Re(v), triv)^{\infty}$ with the property that $f_0 |_K \equiv 1$. Then $f_0(u) \in L^1(U)$ if $v \in ({\f a}^*_{\mathbb C})^{-}$ by Lemma \ref{l1noncompact}. In addition, $f_0(u) >0$ for all $u \in U$.
\begin{thm}\label{supbd}
Let $v \in ({\f a}^*_{\mathbb C})^{-}$ and $f \in \pi(v, \sigma)^{\infty}$. Then $|Wh_{f, m}(a k)| \leq  C_{ m} \| f |_K\|_{sup} a^{\rho+v}$.
\end{thm}
Proof:  Observe that 
$$Wh_{f,m}(ak)= \langle \pi(v, \sigma)(k)f, \pi(-v, \sigma)(a^{-1}) \psi_m \rangle= a^{\rho+v} \langle \pi(v, \sigma)(k)f,  \psi_{Ad(a)m} \rangle.$$
Expressed in the noncompact model
$$Wh_{f,m}(ak)=\int_{U} \pi(v, \sigma)(k)f(u) \psi_{Ad(a)m }(u) d u.$$
Every $|\pi(v, \sigma)(k)f(u)|$ is bounded by  $\| f |_K\|_{sup} f_0(u)$ even though $f_0(u)$ is the spherical vector of $\pi(\Re(v), triv)$. Our theorem then follows from Lemma \ref{l1noncompact}. $\Box$ \\
\\ 
{\bf Proof of Theorem C}: The proof is similar to the proof of Theorem A. Let
us consider the element in the universal enveloping algebra $U(\f g)$, $X^{l}= X_{1, 2}^{l_1} X_{2,3}^{l_2} \ldots X_{n-1,n}^{l_{n-1}}$. Then
$$Wh_{\pi(v, \sigma)(X^{l})f,m}(a)= \langle \pi(v,\sigma)(X^{l})f, \pi(-v, \sigma)(a^{-1}) \psi_m \rangle.$$ 
By Theorem \ref{supbd} and  Lemm \ref{diff},
$$\prod_{i=1}^{n-1} 
(\frac{2 \pi m_i a_i}{a_{i+1}})^{l_i} | Wh_{f,m}(a) |$$ is bounded by
$$\|\pi(v, \sigma)(X^{l})f|_K \|_{sup} C_m a^{\rho+v}.$$
Similar statement is true for every $\pi(v, \sigma)(k)f$:
$$\prod_{i=1}^{n-1} 
(\frac{2 \pi m_i a_i}{a_{i+1}})^{l_i} |Wh_{\pi(v, \sigma)(k)f, m}(a)| \leq \|\pi(v, \sigma)(X^{l})\pi(v, \sigma)(k)f|_K \|_{\sup} C_m a^{\rho+v}.$$
Hence 
$$\prod_{i=1}^{n-1} 
(\frac{2 \pi m_i a_i}{a_{i+1}})^{l_i} |Wh_{f, m}(ak)| \leq C_m a^{\rho+v}
\max_{k_0 \in K} (\|\pi(v, \sigma)(X^{l})\pi(v, \sigma)(k_0)f|_K \|_{\sup})$$
$$=C_m a^{\rho+v}
\max_{k_0 \in K} (\|\pi(v, \sigma)(Ad(k_0^{-1})(X^l))f|_K \|_{\sup}).$$
$\Box$ \\
\\

\subsection{Some Remarks on the general case}
Theorem C can be generalized to cover the case $v$ in the boundary of the negative Weyl chamber.  But the exact same statement will not hold. Instead, for $v \in \partial (\f a_{\mathbb C}^{-})$, we have
$$ \prod_{i=1}^{n-1} 
(\frac{a_i}{a_{i+1}})^{l_i} Wh_f(g) \leq C a^{\rho+v-\epsilon},$$
where $\epsilon \in \f (a^*_{\mathbb C})^+ $ and $\sum_{i=1}^n \epsilon_i=0$.
Indeed, we only require that $\epsilon_i>\epsilon_{i+1}$ if $\Re(v_i)=\Re(v_{i+1})$. If $\Re(v_i)< \Re(v_{i+1})$, $\epsilon_i=\epsilon_{i+1}$ will be allowed. Obviously, these bounds will not be the best bounds. For many applications, these bounds should be adequate. We shall not pursue this here. \\
\\
We are still left with induced representations from $GL(1)$ and $GL(2)$-factors.
By the theorem of Kostant, representations with Whittaker model must have the largest Gelfand-Kirillov dimension (\cite{ko}). By Vogan's classification of unitary dual of $GL(n)$ (\cite{vo}), irreducible unitary generic representations are of the form
$$\pi \cong \Ind_{MAN}^G \mathbb \sigma \otimes \chi  $$
with $$MA= \prod_{i=1}^{r_1} GL(2) \prod_{i=1}^{r_2} GL(1), \qquad (2r_1+r_2=n),$$
$$\sigma=\prod_{i=1}^{r_1} D(d_i+v_i, -d_i+v_i),$$
$$ \chi=\prod_{i=1}^{r_2}  \mathbb C_{\chi_i \epsilon_i}.$$
Here $d_i \in \frac{\mathbb Z^+}{2}$, $v_i \in i \mathbb R$, $\chi_i \in \mathbb C$, $\epsilon_i=\pm 1$. The parameter $\chi$ satisfies the condition that the principal series of $GL(r_2)$ induced from $\chi$ is unitary. $D(d_i+v_i, - d_i+v_i)$ is the discrete series with Harish-Chandra parameter $2 d_i$ and central character $\| \det \|^{v_i}$. The bounds for the Whittaker functions of these $\pi$ is more difficult. At the minimum, We can embedded
$\pi$ as a subquotient of a
certain principal series $\pi(v, \sigma)$ such that $v \in cl(\f a_{\mathbb C})^{-}$. We will be able to obtain bounds on $Wh_{f}$ based on the bounds from the principal series $\pi(v, \sigma)$. These bounds can be far from optimal. The correct way to write down the bounds is to use leading exponents of $\pi(v, \sigma)$ (\cite{knapp}).
When $\sigma$ is a unitary representation of $GL(1)$ factors, the leading exponents are simply $a^{v+\rho}$. For the general case, the leading exponents are more complex.  Hypothetically, the Whittaker functions in this general case should at least share the same kind of bounds as the spherical principal series $\pi(0, triv)$. This should be adequate for applications in automorphic forms.

\commentout{
\subsection{Asymptotes of Whittaker Functions over Siegel set}
 We shall merely point out a functional analytic proof by using a Lemma of Harish-Chandra

Its irreducible unitary representations are considerably simpler, and can be worked out using Mackey's theory. See for example \cite{sahi}  and the reference therein. In addition, by the proof of Kirillov conjecture, any irreducible unitary representation of $G$ will be an irreducible unitary representation of $P_n$  }

\section{Perturbation of Group action and Unitary structure}
Let us get back to the setting of Section 2. Let $(\pi, \mc H)$ be an irreducible unitary representation of $G$ with a Whittaker model and $f \in \mc H_K$. We would like to give a proof of Theorem C and show that
$$Wh_{f}|_{GL(n-1)} \in L^2(U_{n-1} \backslash GL(n-1), |\det|^s) \qquad (\forall \,\, s >0).$$
Proof of Theorem C: By Theorem A, $\forall \,\, t \in \mathbb N$,
$$|\det|^{t} Wh_f|_{GL(n-1)}=\prod_{i=1}^{n-1} (a_i)^{t} Wh_f|_{GL(n-1)} \in 
L^2(\mathbb C_{\mathbf 1} \times_{U_{n-1}} GL(n-1)).$$
Observe that for any $s>0$,
$$\int_{U_{n-1} \backslash GL(n-1)} |Wh_f(g)|^2 |\det(g)|^s d [g] \leq 
\int_{U_{n-1} \backslash GL(n-1), |\det g| < 1} |Wh_f(g)|^2  d [g]$$
$$+\int_{U_{n-1} \backslash GL(n-1), |\det g|| \geq 1} |Wh_f(g)|^2 |\det(g)|^{2 t} d [g] < \infty,$$
where $t$ is an integer greater than $\frac{s}{2}$. $\Box$ \\
\\
Now let us perturb the group action of $\pi$. We define
an action of $G$ on $\pi^{\infty}$: $\forall g \in G$,
$$\pi^{s}(g)=|\det g |^s \pi(g).$$
It is easy to check that $(\pi^s, \pi^{\infty})$ is a  group representation of
$G$. It can never be endowed with a pre-Hilbert structure to make $\pi^s$ unitary. But the following theorem says that $\pi^s|_{P_n}$ can be made into a unitary representation. Indeed, we can perturb the unitary Whittaker model to obtain a unitary structure for $\pi^s|_{P_n}$.
\begin{thm} Let $s \geq 0$
The map
$$Wh: f \in \pi^{\infty} \rightarrow Wh_f(g)|_{P_n} \in L^2(\mathbb C_{\mathbf 1} \times_{U} P_n, |\det |^{2s} d [g])$$
yields a unitary structure of $\pi^s|_{P_n}$ . Here $d [g]= a^{-2 \rho_{n-1}} \frac{ d a}{a} d k$ is the right $P_n$-invariant measure on $ U \backslash P_n$.
\end{thm}
Proof: For any $h \in P_n$ and $f \in \pi^{\infty}$, we compute
\begin{equation}
\begin{split}
& \| Wh_{\pi^s(h)f}(p) \|^2_{L^2(\mathbb C_{\mathbf 1} \times_{U} P_n, |\det |^{2s} d [g])} \\
= & \| Wh_{|\det h|^s \pi(h) f}(P)\|^2_{L^2(\mathbb C_{\mathbf 1} \times_{U} P_n, |\det |^{2s} d [g])} \\
= & |\det h |^{2s} \| R(h) Wh_{f}(p) \|^2_{L^2(\mathbb C_{\mathbf 1} \times_{U} P_n, |\det |^{2s} d [g])} \\
= & \int |\det h|^{2s} | Wh_f(p h) |^2 |\det p|^{2s} d [p] \\
= & \int |\det p|^{2s} |Wh_f(p)|^2 d [p] \\
= & \|Wh_{f}(p) \|^2_{L^2(\mathbb C_{\mathbf 1} \times_{U} P_n, |\det |^{2s} d [g])}
\end{split}
\end{equation}
Here $R(h)$ stands for the right regular action of $h$. 
Hence $\pi^{s}(h)$ preserves the Hilbert norm of $Wh_f|_{P_n}$ in $L^2(\mathbb C_{\mathbf 1} \times_{U} P_n, |\det |^{2s} d [g])$. In fact, $\pi^s(h)$ is simply
$|\det h|^s R(h) $ in the Kirillov-Whittaker model. It is a unitary operator 
of $L^2(\mathbb C_{\mathbf 1} \times_{U} P_n, |\det |^{2s} d [g])$
. $\Box$. \\
\\
Write $\mc H^s=L^2(\mathbb C_{\mathbf 1} \times_{U} P_n, |\det |^{2s} d [g])$. Then $(\pi^s|_{P_n}, \mc H^s)$ is a perturbation of the unitary representation $\pi|_{P_n}$.
\commentout{ The following corollary may also be useful.
\begin{cor} Let $s \geq 0$. For every $f \in \pi^{\infty}$,
$$Wh^s_{f}(p)=\langle \pi^s(p) f, \psi_{\mathbf 1} \rangle \in L^2(\mathbb C_{\mathbf 1} \times_{U} P_{n}, d [p]).$$
\end{cor}
}
In a subsequent paper (\cite{he3}), we shall give some similar results for cusp forms. The philosophy behind this investigation is that many analytic properties for Whittaker functions also hold for cusp forms.

\end{document}